\documentclass[12pt,intlim,righttag]{article}
\usepackage{amsmath,amsthm}
\usepackage{amssymb}
\usepackage{amsfonts}

\newcommand{\la}{\langle}
\newcommand{\ra}{\rangle}

\newcommand\beq{\begin{equation}}
\newcommand\eeq{\end{equation}}

\theoremstyle{Theorem}

\theoremstyle{corollary}

\theoremstyle{remark}

\theoremstyle{definition}

\begin{document}
\title{An extension of the Jacod's condition}

\author{B. Chikvinidze }

\date{~}
\maketitle

\begin{center}

Institute of Cybernetics of Georgian Technical University, 

\
\\       
  Georgian-American University, Business school, 8 M. Aleksidze Srt., 
\\
Tbilisi 0160, Georgia 

\
\\
E-mail: beso.chiqvinidze@gmail.com 

\end{center}

\begin{abstract}


{\bf Abstract} We generalize Jacod's condition and introduce a new type sufficient condition for the uniform integrability of the general stochastic exponential.

\bigskip


\


\end{abstract}

\
\\
{\bf 1. Introduction. }Let us introduce a basic probability space $\big( \Omega , \mathcal{F}, P \big)$ and right continuous filtration 
$(\mathcal{F}_t)_{0\leq t \leq \infty}$ satisfying usual conditions. Let $\mathcal{F}_\infty $ be the smallest $\sigma-$ Algebra containing all $\mathcal{F}_t $ for $t<\infty$ and let $M = (M_t)_{t \geq 0}$ be a local martingale on the stochastic interval $[[0;T]]$ where $T$ is a stopping time. Denote 
by $\bigtriangleup M_t = M_t - M_{t-}$ jumps of $M$ and
by $\mathcal E(M)$ the stochastic exponential of the local martingale $M$:

$$
\mathcal E_{t}(M)= \exp\big\{ M_t-\frac{1}{2}\la M^c\ra_t \big\} \prod_{0<s \leq t}(1+\bigtriangleup M_s)e^{-\bigtriangleup M_s}.
$$
where $M^c$ denotes continuous martingale part of $M$. It is well known that $\mathcal E_{t}(M) = 1 + \int^t_0 \mathcal E_{s-}(M) d M_s$, so it is clear that for local martingale $M$ the associated stochastic exponential $\mathcal E_{t}(M)$ is a local martingale, but not necessarily a true martingale. To know whether $\mathcal E(M)$ is a true martingale is important for many applications, e. g., when Girsanov's transformation is applied to perform a change of measure. Throughout of this paper we assume that 
$\bigtriangleup M_t > -1$ which implies that $\mathcal E_t(M)\geq 0$. So 
$\mathcal E(M)$ will be a supermartingale and the martingale property of $\mathcal E(M)$ is equivalent to $E\mathcal{E}_T(M)=1$.  

\
In case of continuous local martingale $M$ uniform integrability of 
$\mathcal E(M)$ was studied by many authors. Through this paper we are concerned on right continuous exponential martingales. In 1978 Memin and Shiryaev \cite{1} proved that if the elements of the triplet of predictable characteristics of local martingale $M$ are bounded then stochastic exponential $\mathcal E(M)$ is a true martingale. 
Then Lepingle and Memin \cite{2} proved this assertion when the compensator of the process $\la M^c \ra_t + \sum_{s\leq t}(\bigtriangleup M_s)^21_{\{|\bigtriangleup M_s|\leq 1\}} + \sum_{s\leq t} \bigtriangleup M_s 1_{\{|\bigtriangleup M_s| > 1\}}$ is bounded. After that they generalized their result and showed that $Ee^{B_\infty} < \infty$, where $B_t$ is a compensator of the process $A_t = \frac{1}{2}\la M^c \ra_t + \sum_{s\leq t} \{(1+\bigtriangleup M_s)\ln (1+\bigtriangleup M_s) - \bigtriangleup M_s \}$, is sufficient for the uniform integrability of $\mathcal E(M)$. Then J. Jacod \cite{3} introduced sufficient condition in terms of $M^c$ and $\bigtriangleup M$: 
$E\exp \{\frac{1}{2}\la M^c \ra_\infty + \sum_{s\leq \infty }
(\ln (1+\bigtriangleup M_s) - \frac{\bigtriangleup M_s}{1+\bigtriangleup M_s})\} < \infty $. In 2008 P. Protter and K. Shimbo \cite{4} gave sufficient condition in terms of $\la M^c \ra$ and $\la M^d \ra$: $E\exp \{\frac{1}{2}\la M^c \ra_\infty + \la M^d \ra_\infty \} < \infty $. Besides, they showed that the constant $1$ of $\la M^d \ra$ can not be improved (so as the constant $\frac{1}{2}$ at $\la M^c \ra $, see Novikov \cite{6}). In this paper we introduce new type sufficient condition using predictable process $a_s$. Using similar type condition in \cite{5} Chikvinidze obtained necessary and sufficient condition for the uniform integrability of stochastic exponential in case of continuous exponential martingales.

\
\\
Now we formulate the main result of this paper:

\
\\
{\bf Theorem 1}\; Let $M$ be a local martingale with 
$\bigtriangleup M_t > -1$. If there exists some predictable, $M$-integrable process $\; a_s \in [0;1]$ and a constant $\; \varepsilon \;$ with 
$0 < \varepsilon < 1$ such that

$$D =\sup_{0\leq \tau \leq T} E\exp \bigg\{ {\int ^\tau_0 a_s dM_s + \int ^\tau_0  \Big( \frac{1}{2} - a_s \Big) d\langle M^c \rangle_s} + \varepsilon \int^{\tau}_0 1_{\{1-a_s<\varepsilon\}}d \la M^c \ra _s$$

\begin{equation}
 + \sum_{0<s\leq \tau} \Big( \ln(1+\bigtriangleup M_s) - \frac{\bigtriangleup M_s}{1 +\bigtriangleup M_s} + 
\ln (1+a_s\bigtriangleup M_s) - a_s\bigtriangleup M_s \Big) \bigg\} < \infty  
\end{equation}
\\
where $sup$ is taken over all stopping times, then the stochastic exponential $\mathcal E (M)$, defined on the stochastic interval $[[0;T]]$,  is a uniformly integrable martingale.

\
\\
{\bf Remark } If we take $a_s\equiv 0$ then condition $(1)$ from Theorem 1 will turn to Jacod's \cite{3} condition. This means that condition $(1)$ is more general than Jacod's \cite{3} condition.   

\
\\
In section $3$ we construct three counterexamples such that for them the Jacod's \cite{3} condition fails but conditions of Theorem 1 is satisfied for $a_s\equiv 1$ in first example and for $a_s\equiv a\in (0; 1]$ in the second example. Third counterexample shows us the advantage of using predictable processes $a_s$ rather than constant $a$. More precisely, we construct a local martingale such that for any constant 
$a\in [0;1]$ the condition $(1)$ of Theorem 1 fails (therefore Jacod's \cite{3} condition also is not satisfied), but there exists a predictable process $a_s\in [0;1]$ such that conditions of Theorem 1 is satisfied. In all counterexamples the constructed local martingales are purely discontinuous with one jump in first and second cases and with two jumps in the third case. 

\
\\
Finally, in the Appendix, we prove several auxiliary Lemmas used in the proof of the Theorem 1.

\
\\
{\bf 2. Proof of the main result.} Let first prove the following assertion which will be essentially used in the proof of the Theorem 1. 

\
\\
{\bf Proposition 1} \; Under conditions of the Theorem 1 $E\mathcal E_{T}(\int adM)=1$.
\
\\
{\it Proof:} \; It is clear that $a_s\in [0;1]$ and 
$\bigtriangleup M_s > -1$ implies $1+a_s\bigtriangleup M_s> -1$.
So according to Lemma 1 from appendix it is sufficient to show that 
$$\sup_{0\leq \tau \leq T} \bigg[ E\mathcal E_{\tau}\Big(\int adM \Big)\Big[ \frac{1}{2}\int^{\tau}_0 a_s^2\langle M^c \rangle_{s} $$ 
$$ + \sum_{0<s\leq \tau} \Big( \ln(1+a_s\bigtriangleup M_s) - \frac{a_s\bigtriangleup M_s}{1 +a_s \bigtriangleup M_s} \Big) \Big] \bigg]<\infty.$$
For this we will show that the following inequality holds true for some constant $G>0$:

$$\mathcal E_{t}\Big(\int adM \Big)\Big[ \frac{1}{2}\int^{t}_0 a_s^2\langle M^c \rangle_{s} + \sum_{0<s\leq t} \Big( \ln(1+a_s\bigtriangleup M_s) - \frac{a_s\bigtriangleup M_s}{1 +a_s \bigtriangleup M_s} \Big) \Big] $$

$$ \leq G\exp \bigg\{ {\int ^t_0 a_s dM_s + \int^t_0  \Big( \frac{1}{2} - a_s \Big) d\langle M^c \rangle_s} + \varepsilon \int^{t}_0 1_{\{1-a_s<\varepsilon\}}d \la M^c \ra _s $$

$$ + \sum_{0<s\leq t} \Big[ \ln(1+\bigtriangleup M_s) - \frac{\bigtriangleup M_s}{1 +\bigtriangleup M_s} + \ln (1+a_s\bigtriangleup M_s) -a_s\bigtriangleup M_s \Big] \bigg\}.$$

The last inequality is equivalent to the following:

$$ \frac{1}{2}\int^t_0 a_s^2 d\langle M^c \rangle_s + \sum_{0<s\leq t} \Big[ \ln(1+a_s\bigtriangleup M_s) - \frac{a_s\bigtriangleup M_s}{1 +a_s \bigtriangleup M_s} \Big] $$

$$ \leq G\exp \Big\{ \frac{1}{2}\int^t_0 (a_s -1)^2 d\langle M^c \rangle_s + \varepsilon \int^{t}_0 1_{\{1-a_s<\varepsilon\}}d \la M^c \ra _s$$
$$ +\sum_{0<s\leq t} \Big( \ln(1+\bigtriangleup M_s) - \frac{\bigtriangleup M_s}{1 +\bigtriangleup M_s} \Big) \Big\}.$$

Now taking logarithms on the both sides of the inequality we will have: 
$$\ln\bigg( \frac{1}{2}\int^t_0 a_s^2 d\langle M^c \rangle_s + \sum_{0<s\leq t} \Big[ \ln(1+a_s\bigtriangleup M_s) - \frac{a_s\bigtriangleup M_s}{1 +a_s \bigtriangleup M_s} \Big] \bigg) $$
$$\leq \ln G + \frac{1}{2} \int^t_0 (a_s - 1)^2 d\langle M^c \rangle_s + \varepsilon \int^{t}_0 1_{\{1-a_s<\varepsilon\}}d \la M^c \ra _s $$

\begin{equation}
+ \sum_{0<s\leq t} \Big[ \ln(1+\bigtriangleup M_s) - 
\frac{\bigtriangleup M_s}{1 + \bigtriangleup M_s} \Big].
\end{equation}
\\
So, to complete the proof of Proposition 1, it is sufficient to show the validity of the inequality $(2)$. To prove $(2)$ let us first use inequality 
$\ln x \leq \varepsilon^2 x + \ln G$ for some constant $G>0$: 

$$\ln\bigg( \frac{1}{2}\int^t_0 a_s^2 d\langle M^c \rangle_s + \sum_{0<s\leq t} \Big[ \ln(1+a_s\bigtriangleup M_s) - \frac{a_s\bigtriangleup M_s}{1 +a_s \bigtriangleup M_s} \Big] \bigg)$$

\begin{equation}
\leq \ln G + \frac{\varepsilon^2}{2}\int^t_0 a_s^2 d\langle M^c \rangle_s + \varepsilon^2 \sum_{0<s\leq t} \Big[ \ln(1+a_s\bigtriangleup M_s) - 
\frac{a_s\bigtriangleup M_s}{1 + a_s\bigtriangleup M_s} \Big]
\end{equation}
\\
According to Lemma 2 the following inequality holds true:
$$\frac{\varepsilon^2}{2}\int^t_0 a_s^2 d\langle M^c \rangle_s \leq \frac{1}{2}\int^t_0 (a_s-1)^2 d\langle M^c \rangle_s + \varepsilon \int^{t}_0 1_{\{1-a_s<\varepsilon\}}d \la M^c \ra _s.$$
With this because $\ln(1 + a_s\bigtriangleup M_s)-\frac{a_s\bigtriangleup M_s}{1+a_s\bigtriangleup M_s} \geq 0$ we can write:
$$ (3) \leq \ln G + \frac{1}{2}\int^t_0 (a_s-1)^2 d\langle M^c \rangle_s + \varepsilon \int^{t}_0 1_{\{1-a_s<\varepsilon\}}d \la M^c \ra _s $$
$$+ \sum_{0<s\leq t} \Big[ \ln(1+a_s\bigtriangleup M_s) - 
\frac{a_s\bigtriangleup M_s}{1 + a_s\bigtriangleup M_s} \Big].$$
Now inequality $(2)$ follows from the last inequality since for $a_s \in [0;1]$

$$ \ln(1+a_s\bigtriangleup M_s) - 
\frac{a_s\bigtriangleup M_s}{1 + a_s\bigtriangleup M_s}
\leq \ln(1 + \bigtriangleup M_s) - 
\frac{\bigtriangleup M_s}{1 + \bigtriangleup M_s}. $$

\qed

\
\\
Now we are ready to prove the main theorem:
\
\\
{\it Proof of the Theorem 1: }  According to the above mentioned proposition 
$E\mathcal E_{T}(\int adM)=1$. So we can define new probability measure: $dP^a=\mathcal E_{T}(\int adM)dP$. 
It is obvious that $N_t = \int^t_0 (1-a_s)dM_s$ is a $P$ local martingale, so according to the Girsanov's theorem 
$$\tilde{N}_t = N_t - \int^t_0 \frac{1}{1+a_s\bigtriangleup M_s} d \bigg[ \int adM , \int (1-a)dM \bigg]_s $$
$$ = \int^t_0 (1-a_s)dM_s - \int^t_0 \frac{a_s(1-a_s)}{1+a_s\bigtriangleup M_s} d [M]_s$$
will be a $P^a$ local martingale. It is easy to check that 
$\bigtriangleup N_t > -1$ and $\bigtriangleup \tilde{N}_t =
\frac{(1-a_t)\bigtriangleup M_t}{1+a_t\bigtriangleup M_t}> -1 $. Now we show that for local $P^a-$martingale $\tilde{N}$ Jacod's \cite{3} condition is satisfied:
$$E^{P^a} \exp \Big\{\frac{1}{2}\langle \tilde{N}^c \rangle_T + \sum_{0<s\leq T} \big( \ln (1+\bigtriangleup \tilde{N}_s) - \frac{\bigtriangleup \tilde{N}_s}{1+\bigtriangleup \tilde{N}_s}\big) \Big\} $$

$$ = E\mathcal{E}_T \Big( \int adM \Big) \exp \bigg\{ \frac{1}{2}\int^T_0 (1 - a_s)^2 d\langle M^c \rangle_s
$$

\begin{equation}
+ \sum_{0<s\leq T}\Big[ \ln \Big(1+\frac{\bigtriangleup M_s -a_s\bigtriangleup M_s}{1+a_s\bigtriangleup M_s}\Big) - \frac{\frac{(1-a_s)\bigtriangleup M_s}{1+a_s\bigtriangleup M_s}}{1+\frac{(1-a_s)\bigtriangleup M_s}{1+a_s\bigtriangleup M_s}} \Big] \bigg\}
\end{equation}
\\
Here we used equality $\bigtriangleup \tilde{N}_s = \frac{(1-a_s)\bigtriangleup M_s}{1+a_s\bigtriangleup M_s}$. Simplifying $(4)$ we will have:

$$(4) = E\mathcal{E}_T \Big( \int adM \Big)\exp \bigg\{ \frac{1}{2}\int^T_0 (1 - a_s)^2 d\langle M^c \rangle_s $$

\begin{equation}
+ \sum _{0<s\leq T} \Big[ \ln (1+\bigtriangleup M_s) - 
\ln ( 1 + a_s\bigtriangleup M_s ) -
\frac{(1-a_s)\bigtriangleup M_s}{1 + \bigtriangleup M_s} \Big] \bigg\}.
\end{equation}
\\
Now according to Lemma 3 we have an inequality
$$
-\ln (1+a_s\bigtriangleup M_s) - \frac{(1-a_s)\bigtriangleup M_s}{1+\bigtriangleup M_s} \leq -\frac{\bigtriangleup M_s}{1+\bigtriangleup M_s}
$$
which implies the following:
$$
(5) \leq E\mathcal{E}_T \Big( \int adM \Big)\exp \bigg\{ \frac{1}{2}\int^T_0 (1 - a_s)^2 d\langle M^c \rangle_s
$$

$$
+ \sum _{0<s\leq T} \Big[ \ln (1+\bigtriangleup M_s)  
-\frac{\bigtriangleup M_s}{1 + \bigtriangleup M_s} \Big] \bigg\}\leq D < \infty.
$$
\\
This means that $\mathcal{E}(\tilde{N})$ is a $P^a$ uniformly integrable martingale on the stochastic interval $[[0;T]]$. So $E^{P^a}\mathcal{E}_T(\tilde{N})=1$. Now let us show that this implies $E\mathcal{E}_T(M)=1$.

$$E^{P^a}\mathcal{E}_T(\tilde{N})= E\bigg[\exp \Big\{ \int^T_0 a_s dM_s - \frac{1}{2}\int^T_0 a_s^2 d\langle M^c \rangle_s \Big\} 
\prod _{0<s\leq T}(1+a_s\bigtriangleup M_s)e^{-a_s\bigtriangleup M_s} $$

$$\times \exp \Big\{ \int^T_0 (1 - a_s) dM_s - \int^T_0 \frac{a_s(1-a_s)}{1+a_s\bigtriangleup M_s} d[M]_s - \frac{1}{2}\int^T_0 (1-a_s)^2 d\langle M^c \rangle_s 
\Big\} $$

$$\times \prod _{0<s\leq T}\Big( 1+\frac{(1-a_s)\bigtriangleup M_s}{1+a_s\bigtriangleup M_s} \Big)e^{-\frac{(1-a_s)\bigtriangleup M_s}{1+a_s\bigtriangleup M_s}} \bigg]$$

$$= E \bigg[ \exp \Big\{ M_T -\frac{1}{2}\int^T_0 a_s^2 d\langle M^c \rangle_s - \int^T_0 \frac{a_s(1-a_s)}{1+a_s\bigtriangleup M_s} d\langle M^c \rangle_s $$

$$ - \sum_{0<s\leq T} \frac{a_s(1-a_s)}{1+a_s\bigtriangleup M_s}(\bigtriangleup M_s)^2 - \frac{1}{2}\int^T_0 (a_s-1)^2 d\langle M^c \rangle_s  \Big\} $$

\begin{equation}
\times \prod _{0<s\leq T} (1+a_s\bigtriangleup M_s)\Big( 1+\frac{(1-a_s)\bigtriangleup M_s}{1+a_s\bigtriangleup M_s} \Big)e^{-a_s\bigtriangleup M_s -\frac{(1-a_s)\bigtriangleup M_s}{1+a_s\bigtriangleup M_s}}\bigg].
\end{equation}
\\
It is clear that $\int^T_0 \frac{a_s(1-a_s)}{1+a_s\bigtriangleup M_s} d\langle M^c \rangle_s = \int^T_0 a_s(1-a_s) d\langle M^c \rangle_s$, because $\la M^c \ra$ is continuous. Using this we will have from $(6)$: 

$$(6)= E\bigg[\exp \Big\{ M_T - \frac{1}{2} \int^T_0 a_s^2 \langle M^c \rangle_s - \int^T_0 (a_s - a_s^2) \langle M^c \rangle_s $$

$$
-\frac{1}{2} \int^T_0 (a_s - 1)^2 \langle M^c \rangle_s \Big\} \times \prod _{0<s\leq T} e^{-\frac{a_s(1-a_s)}{1+a_s\bigtriangleup M_s} (\bigtriangleup M_s)^2}$$

$$\times \prod _{0<s\leq T}\big(1+a_s\bigtriangleup M_s + (1-a_s)\bigtriangleup M_s \big) e^{\frac{-a_s\bigtriangleup M_s - a_s^2 (\bigtriangleup M_s)^2 - \bigtriangleup M_s + a_s\bigtriangleup M_s}{1+a_s\bigtriangleup M_s}} \bigg]    $$

$$=E \bigg[ \exp \big\{ M_T -\frac{1}{2}\la M^c \ra_T \big\}\prod _{0<s\leq T} (1+\bigtriangleup M_s)
e^{\frac{-\bigtriangleup M_s - a_s^2 (\bigtriangleup M_s)^2 - a_s (\bigtriangleup M_s)^2 + a_s^2 (\bigtriangleup M_s)^2}{1+a_s\bigtriangleup M_s}} \bigg]$$

$$=E \bigg[ \exp \big\{ M_T -\frac{1}{2}\la M^c \ra_T \big\}\prod _{0<s\leq T} (1+\bigtriangleup M_s)
e^{-\bigtriangleup M_s } \bigg]=E\mathcal{E}_T(M).$$
\\
So we obtained that $E\mathcal{E}_T(M)=E^{P^a}\mathcal{E}_T(\tilde{N})=1$. This means that 
$\mathcal{E}(M)$ is a uniformly integrable martingale on the stochastic interval $[[0;T]]$.
\qed

\
\\
{\bf 3. Counterexamples.}
\\
{\bf I.} \; In first part we construct discrete-time martingale for which Jacod's \cite{3} condition fails but conditions of Theorem 1 holds true for constant $a_s\equiv 1$. Let us consider continuous random variable $\xi$ with the probability density function $f(x)$:

$$
f(x)=
\begin{cases}
{\frac{1}{2(1+x)^2}e^{\frac{x}{1+x}}\;\;\; : \;\;\; -1 < x \leq 0},\\
{\frac{1}{2(1-x)^2}e^{-\frac{x}{1-x}}\;\;\; : \;\;\; 0\leq x<1},\\
{\;\;\;\;\;\;0 \;\;\;\;\;\;\;\;\;\;  :  \;\;\;\;\;\;\;\; otherwise}
\end{cases}
$$
It easy to check that $\int_{-1}^{1}f(x)dx = 1$, $E|\xi|<\infty $ and $E\xi =0$. Define the discrete flow of $\sigma$-Algebras: 
$\mathcal{F}_0 = \{\emptyset ; \Omega \}$ and $\mathcal{F}_k = \sigma (\xi)$ for any $k\geq 1$. Then $M_t = \xi 1_{\{t\geq 1\}}$ will be a discrete time martingale with one jump 
$\bigtriangleup M_1 = \xi$. Now we will show that for $M$ Jacod's \cite{3} condition is not satisfied. It is obvious that $M^c\equiv 0$ and $\bigtriangleup M_k = \xi 1_{\{k=1\}}$ so we will have:

$$
E\exp\Big\{\sum _{k}\big[\ln (1+\bigtriangleup M_k) - \frac{\bigtriangleup M_k}{1+\bigtriangleup M_k}\big]\Big\} = E\exp \big\{\ln (1+\xi)-\frac{\xi}{1+\xi}\big\}
$$

$$
=E(1+\xi)e^{-\frac{\xi}{1+\xi}} = \int_{-1}^{1} (1+x)e^{-\frac{x}{1+x}} f(x) dx > \int_{-1}^{0} (1+x)e^{-\frac{x}{1+x}} \frac{1}{2(1+x)^2}e^{\frac{x}{1+x}} dx 
$$

$$
=\frac{1}{2} \int_{-1}^{0} \frac{1}{1+x} dx = 
\frac{1}{2} \ln(1+x) \bigg|^0_{-1} = \infty .
$$
Now let us check condition $(1)$ of Theorem 1 when $a_s\equiv 1$. For any stopping time $\tau$ we will have:

$$
E\exp\Big\{M_{\tau} + \sum _{0<k\leq \tau}\big[\ln (1+\bigtriangleup M_k) - \frac{\bigtriangleup M_k}{1+\bigtriangleup M_k} + \ln(1+\bigtriangleup M_k) - \bigtriangleup M_k \big]\Big\}
$$

$$
\leq E\exp \big\{\ln (1+\xi)-\frac{\xi}{1+\xi} + \ln (1+\xi)\big\} = E(1+\xi)^2 e^{-\frac{\xi}{1+\xi}}
$$

$$
=\frac{1}{2} \int_{-1}^0 (1+x)^2 e^{-\frac{x}{1+x}} \times \frac{1}{(1+x)^2}e^{\frac{x}{1+x}} dx + \frac{1}{2}\int_{0}^1 (1+x)^2 e^{-\frac{x}{1+x}} \times \frac{1}{(1-x)^2}e^{-\frac{x}{1-x}} dx
$$

$$
\leq \frac{1}{2} + 4e^2 \int_0^1 \frac{1}{(1-x)^2} e^{-\frac{2}{1-x}} dx = 2e^2 \big( -e^{-\frac{2}{1-x}} \big) \Big|_0^1 = 
2e^2(0+e^{-2})=2.
$$

\
\\
{\bf II.} \; In second part we construct local martingale for which Jacod's \cite{3} condition fails but conditions of Theorem 1 holds true for any constant $a_s\equiv a\in (0;1]$. Let $N_t$ be a standard Poisson martingale and 
$\tau _1$  be the first jump moment for $N$. Consider the martingale $M_t=\int_0^t e^s dN_{s\wedge \tau _1}$. It is clear that $M$ has jump only in $\tau _1$ and 
$\bigtriangleup M_{\tau _1} = e^{\tau _1}$. First let us check that Jacod's \cite{3} condition fails for $M$. Because $M^c\equiv 0$ we will have:

$$
E\exp\Big(\sum _{0<s\leq \tau _1}\big[\ln (1+\bigtriangleup M_s) - \frac{\bigtriangleup M_s}{1+\bigtriangleup M_s}\big]\Big) = E\exp \Big(\ln (1+e^{\tau _1})-\frac{e^{\tau _1}}{1+e^{\tau _1}}\Big)
$$

$$
\geq e^{-1}E(1+e^{\tau _1})=e^{-1}+e^{-1}Ee^{\tau _1}=
e^{-1}+e^{-1}\int_0^\infty e^t\times e^{-t} dt = \infty .
$$
Now we show that martingale $M$ satisfies condition $(1)$ of Theorem 1 for $a_s\equiv a$ where $a$ is any constant from the interval $(0;1]$. For any stopping time $\tau $ we will have:

$$
E\exp \Big\{ aM_{\tau} + \sum _{0<s\leq \tau}\big[\ln (1+\bigtriangleup M_s) - \frac{\bigtriangleup M_s}{1+\bigtriangleup M_s} + \ln(1+a\bigtriangleup M_s) - a\bigtriangleup M_s \big] \Big\}
$$ 

$$
\leq E\exp \Big\{ aM_{\tau _1} + \ln (1+\bigtriangleup M_{\tau _1}) - \frac{\bigtriangleup M_{\tau _1}}{1+\bigtriangleup M_{\tau _1}} + \ln(1+a\bigtriangleup M_{\tau _1}) - a\bigtriangleup M_{\tau _1} \Big\}
$$

$$
=E \exp \Big\{ ae^{\tau _1} -a\int_0^{\tau _1} e^s ds 
+ \ln (1+e^{\tau _1}) - \frac{e^{\tau _1}}{1+e^{\tau _1}} +
\ln (1+ae^{\tau _1}) -ae^{\tau _1}  \Big\}
$$

\begin{equation}
\leq e^a E \exp \big\{ -ae^{\tau _1} + \ln (1+e^{\tau _1}) +
\ln (1+ae^{\tau _1}) \big\}
\end{equation}
If we use inequality $\ln x \leq \delta x + G$ where 
$\delta \in (0;\frac{a}{1+a})$ and $G>0$ are some constants, we obtain from $(7)$:

$$
(7)\leq e^{a}E \exp \{-ae^{\tau _1} + 2\delta +2G + \delta e^{\tau _1} + 
a\delta e^{\tau _1}\} $$
$$= e^{a+2\delta +2G}
E\exp \{-(a-(1+a)\delta )e^{\tau _1}\} \leq e^{a+2\delta +2G}< \infty .
$$

\
\\
{\bf III.} \; This counterexample considered now shows us advantage of predictable process $a_s\in [0;1]$ rather than constants $a\in [0;1]$. 

\
Let us consider a random variable $\eta $ with the probability density function $g(x)$:

$$
g(x)=
\begin{cases}
{-3x+1\;\;\; : \;\;\;   x \in [-\frac{1}{2};0]},\\
{\frac{1}{4x^3}\;\;\; : \;\;\; x \in [1;\infty )},\\
{0 \;\;\; :   \;\;\;\;\;\; otherwise}
\end{cases}
$$
It easy to check that $\int_{-\infty}^{+\infty} g(x)dx = 1$, $E|\eta |<\infty$, $E\eta =0$ and $E\eta^2 = \infty $. With this let $N_t$ be the standard Poisson martingale which is independent from random variable $\eta$. Denote by $\hat{N}_t = N_t - N_{t\wedge 1}$ standard Poisson martingale starting from $1$ and let 
$\hat{\tau}_1$ be the first jump moment for $\hat{N}_t$. 
\

Define a flow of $\sigma$-Algebras: 
$\mathcal{F}_t = \{\emptyset ; \Omega \}$ if $0\leq t <1$ and 
$\mathcal{F}_t = \sigma (\eta , \hat{N}_s \; : \; s\leq t)$ for $t \geq 1$. Define $M^1_t=\eta 1_{\{t\geq 1\}}$ and 
$M^2_t = \int_1^t e^{s-1}d \hat{N}_{s\wedge \hat{\tau}_1}$. It is obvious that $M^1$ and $M^2$ are local martingales with respect to the filtration $(\mathcal{F}_t)_{t\geq 0}$.
\

In second counterexample we proved that for $M^2$ the Jacod's \cite{3} condition fails, but for every constant $a\in (0;1]$ the conditions of Theorem 1 is satisfied. First we will show that for $M^1$ the Jacod's \cite{3} condition holds true, but for every constant $a\in (0;1]$ the condition $(1)$ of Theorem 1 fails.

\begin{equation}
Ee^{\ln (1+\bigtriangleup M^1_1)-\frac{\bigtriangleup M^1_1}{1+\bigtriangleup M^1_1}}=e^{-1}Ee^{\ln (1+\eta)+\frac{1}{1+\eta}}
\leq eE(1+\eta)=e+eE\eta =e.
\end{equation}

Now taking any constant $a\in (0;1]$, we will have: 
$$
E \exp \Big\{ aM^1_1 + \ln (1+\bigtriangleup M^1_1) - \frac{\bigtriangleup M^1_1}{1+\bigtriangleup M^1_1} + 
\ln (1+a\bigtriangleup M^1_1) -a\bigtriangleup M^1_1 \Big\}
$$

$$
=e^{-1}E\exp \Big\{ a\eta + \ln (1+\eta ) + \frac{1}{1+\eta}
+\ln (1+a\eta ) -a\eta  \Big\}
$$

\begin{equation}
\geq e^{-1}E\exp \{\ln (1+\eta) + \ln (1+a\eta)\}=e^{-1}E(1+\eta)(1+a\eta)=e^{-1}(1+aE\eta^2)=\infty.
\end{equation}

Define local martingale $M$:
$$M_t=M^1_t+M^2_t=\eta 1_{\{t\geq 1\}} + \int_1^t e^{s-1}d \hat{N}_{s\wedge \hat{\tau}_1}.$$
It is obvious that $M$ has two jumps at $1$ and at $\hat{\tau}_1$. Now we will show that for any constant $a\in [0;1]$ condition $(1)$ of Theorem 1 fails. 

$$
E\exp \Big\{ aM_{\hat{\tau}_1} + \sum _{0<s\leq \hat{\tau}_1}\big[\ln (1+\bigtriangleup M_s) - \frac{\bigtriangleup M_s}{1+\bigtriangleup M_s} + \ln(1+a\bigtriangleup M_s) - a\bigtriangleup M_s \big] \Big\}
$$

$$
=E\exp \Big\{ a\eta + \ln (1+\eta) - \frac{\eta}{1+\eta} + \ln(1+a\eta) - a\eta \Big\}
$$

\begin{equation}
\times E \exp \Big\{ ae^{\hat{\tau}_1-1} -a\int_1^{\hat{\tau}_1} e^{s-1} ds 
+ \ln (1+e^{\hat{\tau}_1-1}) - \frac{e^{\hat{\tau}_1-1}}{1+e^{\hat{\tau}_1-1}} +
\ln (1+ae^{\hat{\tau}_1-1}) -ae^{\hat{\tau}_1-1}  \Big\}
\end{equation}
Here we used the independence of $\eta$ and $\hat{\tau}_1$. Simplifying 
$(10)$ we obtain: 

$$
(10)=E\exp \Big\{\ln (1+\eta) - \frac{\eta}{1+\eta} + \ln(1+a\eta) \Big\} 
$$

\begin{equation}
\times e^aE \exp \Big\{ -ae^{\hat{\tau}_1-1} 
+ \ln (1+e^{\hat{\tau}_1-1}) - \frac{e^{\hat{\tau}_1-1}}{1+e^{\hat{\tau}_1-1}} +
\ln (1+ae^{\hat{\tau}_1-1}) \Big\}.
\end{equation}
If $a=0$ as we proved in second counterexample, the second multiplier of 
$(11)$ is infinity, and if $0<a\leq 1$ then according to $(9)$ the first multiplier of $(11)$ will turn to infinity. So for any $a\in [0;1]$ condition $(1)$ of Theorem 1 fails for local martingale $M$.

\
\\
Now consider predictable process $a_s=1_{\{s>1\}}$. It is obvious that $a_s\in [0;1]$. So we left to check condition $(1)$ of Theorem 1. For any stopping time $\tau$ we will have: 

$$
E\exp \Big\{ \int_0^{\tau}a_sdM_s + \sum _{0<s\leq \tau}\big[\ln (1+\bigtriangleup M_s) - \frac{\bigtriangleup M_s}{1+\bigtriangleup M_s} + \ln(1+a_s\bigtriangleup M_s) - 
a_s\bigtriangleup M_s \big] \Big\}
$$
 
$$
\leq E\exp \Big\{ \int_0^1a_sdM_s + \ln (1+\bigtriangleup M_1) - \frac{\bigtriangleup M_1}{1+\bigtriangleup M_1} + \ln(1+a_1\bigtriangleup M_1) - a_1\bigtriangleup M_1 \Big\}
$$

\begin{equation}
\times \exp \Big\{ \int_1^{\hat{\tau}_1}a_sdM_s + \ln (1+\bigtriangleup M_{\hat{\tau}_1}) - \frac{\bigtriangleup M_{\hat{\tau}_1}}{1+\bigtriangleup M_{\hat{\tau}_1}} + \ln(1+a_{\hat{\tau}_1}\bigtriangleup M_{\hat{\tau}_1}) - a_{\hat{\tau}_1}\bigtriangleup M_{\hat{\tau}_1} \Big\}
\end{equation}
If we simplify $(12)$ and use independence of $\eta$ and $\hat{N}_t$ we obtain: 
$$
(12)=E\exp \Big\{ \ln(1+\eta) -\frac{\eta}{1+\eta} \Big\}
$$

\begin{equation}
\times E\exp \Big\{ - e^{\hat{\tau}_1-1} +1 + 
\ln (1+e^{\hat{\tau}_1-1}) - \frac{e^{\hat{\tau}_1-1}}
{1+e^{\hat{\tau}_1-1}} + \ln(1+e^{\hat{\tau}_1-1}) \Big\}.
\end{equation}
It follows from $(8)$ that first part of $(13)$ is finite and as we have mentioned above, in second counterexample we proved that the second part of $(13)$ is also finite. This means that for local martingale $M$ the condition $(1)$ of Theorem 1 fails for any constant $a\in [0;1]$, but for the process $a_s=1_{\{s>1\}}$ - condition $(1)$ is satisfied.

\
\\
{\bf 4. Appendix.}
\\
{\bf Lemma 1} \; Let $M$ be a local martingale with 
$\bigtriangleup M_s > -1$. If 
$$ \sup_{0\leq \tau \leq T} E \bigg[ \mathcal E_{\tau}(M)\Big\{ \frac{1}{2}\langle M^c \rangle_{\tau} + \sum_{0<s\leq \tau} \Big( \ln(1+\bigtriangleup M_s) - \frac{\bigtriangleup M_s}{1 +\bigtriangleup M_s} \Big) \Big\} \bigg]<\infty $$
then $\mathcal E(M)$ is a uniformly integrable martingale on $[[0;T]]$.
\
\\
{\it Proof:} \; 
Using Ito's formula we will have:
$$\mathcal E_{t}(M)\ln \mathcal E_{t}(M) = \int^t_0 (\ln\mathcal E_{s-}(M) + 1)\mathcal E_{s-}(M) d M_s + \frac{1}{2}\int^t_0 \mathcal E_{s-}(M) d \langle M^c \rangle_s  $$
$$ + \sum _{0<s\leq t} \big[ \mathcal E_{s}(M)\ln\mathcal E_{s}(M) - \mathcal E_{s-}(M)\ln\mathcal E_{s-}(M) - (\ln \mathcal E_{s-}(M) + 1)\bigtriangleup \mathcal E_{s}(M) \big].$$
Notice that $\bigtriangleup \mathcal E_{t}(M) = \mathcal E_{t-}(M)\bigtriangleup M_t$ and $\frac{\mathcal E_{t}(M)}{\mathcal E_{t-}(M)} = 1 + \bigtriangleup M_t$. Using these equalities and localization arguments, there exists increasing sequence of stopping times $(\tau _n)_{n\geq 1}$ such that:

$$E\mathcal E_{\tau _n}(M)\ln \mathcal E_{\tau _n}(M) = \frac{1}{2}E\int^{\tau _n}_0 \mathcal E_{s-}(M) d \langle M^c \rangle_s  $$
$$+ E\sum _{0<s\leq \tau_n} \big[ \mathcal E_{s}(M)\ln\mathcal E_{s}(M) - \mathcal E_{s-}(M)\ln\mathcal E_{s-}(M) - \mathcal E_{s}(M)\ln\mathcal E_{s-}(M)  $$
$$ + \mathcal E_{s-}(M)\ln\mathcal E_{s-}(M) - \bigtriangleup \mathcal E_{s}(M) \big]  $$
$$ = \frac{1}{2}E\int^{\tau _n}_0 \mathcal E_{s-}(M) d \langle M^c \rangle_s + E\sum _{0<s\leq \tau_n} \big[ \mathcal E_{s}(M)\ln(1 + \bigtriangleup M_s) - 
\mathcal E_{s-}(M)\bigtriangleup M_s \Big] $$
\begin{equation}
= \frac{1}{2}E\int^{\tau _n}_0 \mathcal E_{s-}(M) d \langle M^c \rangle_s + E\sum _{0<s\leq \tau_n} \mathcal E_{s-}(M) \big[ (1 + \bigtriangleup M_s)\ln(1 + \bigtriangleup M_s) - \bigtriangleup M_s \Big].
\end{equation}
 
On the other hand 
$$\mathcal E_{t}(M)\Big[ \frac{1}{2}\langle M^c \rangle_t + \sum_{0<s\leq t} \Big( \ln(1+\bigtriangleup M_s) - \frac{\bigtriangleup M_s}{1 +\bigtriangleup M_s} \Big) \Big] =local \; martingale$$

$$ + \frac{1}{2}\int^t_0 \mathcal E_{s-}(M) d \langle M^c \rangle_s + \sum_{0<s\leq t} \mathcal E_{s-}(M) \Big( \ln(1+\bigtriangleup M_s) - \frac{\bigtriangleup M_s}{1 +\bigtriangleup M_s} \Big) $$

$$ + \sum_{0<s\leq t} \Big( \ln(1+\bigtriangleup M_s) - \frac{\bigtriangleup M_s}{1 +\bigtriangleup M_s} \Big)\bigtriangleup \mathcal E_{s}(M).$$
\\
Again using localization arguments we obtain: 

$$E\mathcal E_{\tau_n}(M)\Big[ \frac{1}{2}\langle M^c \rangle_{\tau_n} + \sum_{0<s\leq \tau_n} \Big( \ln(1+\bigtriangleup M_s) - \frac{\bigtriangleup M_s}{1 +\bigtriangleup M_s} \Big) \Big] = \frac{1}{2}E\int^{\tau _n}_0 \mathcal E_{s-}(M) d \langle M^c \rangle_s $$

$$+ E\sum_{0 < s \leq \tau_n} \mathcal E_{s-}(M) \Big[ \ln(1+\bigtriangleup M_s) - \frac{\bigtriangleup M_s}{1 +\bigtriangleup M_s} +
\bigtriangleup M_s \times \ln(1 + \bigtriangleup M_s) - \frac{(\bigtriangleup M_s)^2}{1 +\bigtriangleup M_s} \Big] $$
\begin{equation}
= \frac{1}{2}E\int^{\tau _n}_0 \mathcal E_{s-}(M) d \langle M^c \rangle_s + E\sum_{0 < s \leq \tau_n} \mathcal E_{s-}(M) 
\Big[ (1+\bigtriangleup M_s)\ln(1+\bigtriangleup M_s) - \bigtriangleup M_s \Big].
\end{equation}
\\  
So finally we get from $(14)$ and $(15)$ that: 

$$E\mathcal E_{\tau _n}(M)\ln \mathcal E_{\tau _n}(M) = \frac{1}{2}E\int^{\tau _n}_0 \mathcal E_{s-}(M) d \langle M^c \rangle_s $$

$$ + E\sum_{0 < s \leq \tau_n} \mathcal E_{s-}(M) \Big[ (1+\bigtriangleup M_s)\ln(1+\bigtriangleup M_s) - \bigtriangleup M_s \Big] $$

$$=E\mathcal E_{\tau_n}(M)\Big[ \frac{1}{2}\langle M^c \rangle_{\tau_n} + \sum_{0<s\leq \tau_n} \Big( \ln(1+\bigtriangleup M_s) - \frac{\bigtriangleup M_s}{1 +\bigtriangleup M_s} \Big) \Big] < \infty .$$
\\
This means that the family $\big(\mathcal E_{\tau _n}(M)\big)_{n\geq 1}$ is uniformly integrable, so $\mathcal E(M)$ is a uniformly integrable martingale.
\qed 

\
\\
{\bf Lemma 2 }
For any $x \in [0;1]$ and $0<\varepsilon <1$ the following inequality holds:
$$(1-\varepsilon^2)x^2-2x+1+2\varepsilon 1_{\{1-x<\varepsilon \}}\geq 0 .$$
\
\\
{\it Proof:} \; The solutions of the equation 
$(1-\varepsilon^2)x^2-2x+1=0$ are $x_1=\frac{1}{1+\varepsilon}<1$ and $x_2=\frac{1}{1-\varepsilon}>1$. If $x\leq x_1$ then 
$(1-\varepsilon^2)x^2-2x+1\geq 0$ and if $x_1 < x \leq 1$ then $(1-\varepsilon^2)x^2-2x+1\geq -\varepsilon^2$. It is clear that inequality $x_1 < x \leq 1$ implies $1-x<\varepsilon$, so if $x_1 < x \leq 1$ then 
$(1-\varepsilon^2)x^2-2x+1+2\varepsilon 1_{\{1-x<\varepsilon \}}\geq -\varepsilon^2 + 2\varepsilon > 0$ because $0<\varepsilon <1$.
\qed    

\
\\
{\bf Lemma 3 }
If $a_s \in [0;1]$ and $\bigtriangleup M_s > -1$ then  

\begin{equation}
\frac{\bigtriangleup M_s}{1+\bigtriangleup M_s} \leq \ln (1+a_s\bigtriangleup M_s) + \frac{(1-a_s)\bigtriangleup M_s}{1+\bigtriangleup M_s}.
\end{equation}

\
\\
{\it Proof:} \; It is easy to check that $(1+\bigtriangleup M_s)\ln (1+\bigtriangleup M_s) - \bigtriangleup M_s \geq 0$ and 
$\bigtriangleup M_s \ln (1+a_s\bigtriangleup M_s)\geq 
a_s\bigtriangleup M_s \ln (1+a_s\bigtriangleup M_s)$. Combining these two inequalities we obtain:

$$\bigtriangleup M_s \leq (1+a_s\bigtriangleup M_s)\ln (1+a_s\bigtriangleup M_s) - a_s\bigtriangleup M_s + \bigtriangleup M_s$$
$$
\leq (1+\bigtriangleup M_s)\ln(1+a_s\bigtriangleup M_s) + 
(1-a_s)\bigtriangleup M_s. 
$$
If we divide by $1+\bigtriangleup M_s$ the both sides of the last inequality we get $(16)$.

\qed

\end{document}